\documentclass{amsart}
\usepackage{amsmath,amssymb}
\usepackage{yhmath}

\newcommand{\bdis}{\begin{displaymath}}
\newcommand{\edis}{\end{displaymath}}
\newcommand{\be}{\begin{equation}}
\newcommand{\ee}{\end{equation}}
\newcommand{\mbb}{\mathbb}
\newcommand{\mcal}{\mathcal}

\newcommand{\vp}{\varphi}

\newcommand{\zf}{\zeta\left(\frac{1}{2}+it\right)}

\theoremstyle{definition}

\theoremstyle{remark}
\newtheorem{remark}[]{Remark}

\newtheorem*{mydef11}{{\bf Theorem 1}}

\newtheorem*{mydef12}{{\bf Theorem 2}}

\newtheorem*{mydef13}{{\bf Theorem 3}}

\newtheorem*{mydef14}{{\bf Theorem 4}}

\newtheorem*{mydef41}{{\bf Corollary 1}}

\newtheorem*{mydef42}{{\bf Corollary 2}}

\newtheorem*{mydef51}{{\bf Lemma 1}}

\newtheorem*{mydef52}{{\bf Lemma 2}}

\newtheorem*{mydef53}{{\bf Lemma 3}}

\newtheorem*{mydef54}{{\bf Lemma 4}}

\newtheorem*{mydef55}{{\bf Lemma 5}}

\newtheorem*{mydef56}{{\bf Lemma 6}}

\newtheorem*{mydefII1}{\bf Consequence 1}

\numberwithin{equation}{section}

\begin{document}

\title[Jacob's ladders, logarithmic modification \dots]{Jacob's ladders, logarithmic modification of the Hardy-Littlewood integral (1918), Titchmarsh's $\Omega$-theorem (1928) and new point of contact with the Fermat-Wiles theorem}

\author{Jan Moser}

\address{Department of Mathematical Analysis and Numerical Mathematics, Comenius University, Mlynska Dolina M105, 842 48 Bratislava, SLOVAKIA}

\email{jan.mozer@fmph.uniba.sk}

\keywords{Riemann zeta-function}

\begin{abstract}
In this paper we obtain two new points of contact between Jacob's ladders and Fermat-Wiles theorem. They are generated by a logarithmic modification of the Hardy-Littlewood integral. Furthermore, we present a kind of asymptotic laws of conservation for a set of areas connected with above mentioned modification of the Hardy-Littlewood integral. 
\end{abstract}
\maketitle

\section{Introduction} 

\subsection{} 

In our previous papers \cite{6}-- \cite{10} we have studied the Hardy-Littlewood integral 
\be \label{1.1} 
J(T)=\int_0^T\left|\zf\right|^2{\rm d}t. 
\ee  
Namely, we have been interested in some sets of almost linear increments of $J(T)$. Apart from others, we have achieved in this direction a number of the equivalents of the Fermat-Wiles theorem. 

In this paper we present results concerning the following logarithmic modification 
\be \label{1.2} 
J_1(T)=\int_0^T\left\{ \ln t- \left|\zf\right|^2 \right\}{\rm d}t 
\ee 
of the Hardy-Littlewood integral (\ref{1.1}). Namely, we obtain new results about oscillations of the function $\left|\zf\right|^2$ around the curve $y=\ln t$, and also new equivalents of the Fermat-Wiles theorem. 

\subsection{} 

First, we obtain the following: modification (\ref{1.2}) of the Hardy-Littlewood integral generates new point of contact with the Fermat-Wiles theorem. Namely, it is true that the condition 
\be \label{1.3} 
\lim_{\tau\to\infty}\frac{1}{\tau}\int_0^{\frac{x^n+y^n}{z^n}\frac{\tau}{\ln 2\pi-2c}}\left\{ \ln t- \left|\zf\right|^2 \right\}{\rm d}t \not=1
\ee 
for all Fermat's rationals 
\bdis 
\frac{x^n+y^n}{z^n};\ x,y,z,n\in\mbb{N},\ n\geq 3
\edis  
represents logarithmically modified $\zeta$-equivalent of the Fermat-Wiles theorem. 

Next, we give also a second variant of (\ref{1.3}). Namely the condition 
\be \label{1.4} 
\lim_{\tau\to\infty}\frac{1}{\tau}\int_{\frac{x^n+y^n}{z^n}\frac{\tau}{\ln 2\pi-2c}}^{(N+1)\frac{x^n+y^n}{z^n}\frac{\tau}{\ln 2\pi-2c}}\left\{ \ln t- \left|\zf\right|^2 \right\}{\rm d}t\not=N 
\ee 
for all Fermat's rationals 
\bdis 
\frac{x^n+y^n}{z^n} 
\edis 
and for every fixed $N\in\mbb{N}$ (separately) represents the second variant of the logarithmically modified $\zeta$-equivalent of the Fermat-Wiles theorem. 

\subsection{} 

Now, let us remind the Titchmarsh $\Omega$-theorem (see \cite{11}, pp. 78, 79; comp. \cite{12}, p. 172): the inequality 
\be \label{1.5} 
\left|\zf\right|>\exp[\ln^\alpha t]
\ee 
is satisfied for some indefinitely large values of $t$ provided that 
\be \label{1.6} 
\alpha<\frac 12. 
\ee  

Let us define the following sets 
\be \label{1.7} 
\begin{split}
& (0,T]^+=\left\{t:\ \ln t>\left|\zf\right|^2,\ 0<t\leq T\right\}, \\ 
& (0,T]^-=\left\{t:\ \ln t\leq \left|\zf\right|^2,\ 0<t\leq T\right\}, 
\end{split} 
\ee  
for every sufficiently big $T$. 

However, it is true that 
\be \label{1.8} 
\left|\zf\right|\geq e^{\frac 12\ln\ln t},\ t\in (0,T]^-, 
\ee 
and consequently, by the Titchmarsh $\Omega$-theorem, we have that the set $(0,T]^-$ is disconnected set composed by the finite number of intervals. Of course, the set $(0,T]^+$ has similar structure. 

\subsection{} 

Let us remind the formula\footnote{See \cite{2}, comp. \cite{3}.} 
\be \label{1.9} 
\begin{split}
& \int_0^T\left|\zf\right|^2{\rm d}t=\vp_1(T)\ln\{\vp_1(T)\}+(c-\ln 2\pi)\vp_1(T)+ \\ 
& c_0+\mcal{O}\left(\frac{\ln T}{T}\right), 
\end{split}
\ee 
that implies the infinite set of almost exact representations of the Hardy-Littlewood integral (\ref{1.1}). 

Now, the formula (\ref{1.9}) gives the following result in the case (\ref{1.2}): 
\be \label{1.10} 
\begin{split}
& \int_{t\in (0,\overset{r-1}{T}]^+}\left\{ \ln t- \left|\zf\right|^2 \right\}{\rm d}t - \\ 
& \int_{t\in (0,\overset{r-1}{T}]^-}\left\{ \ln t- \left|\zf\right|^2 \right\}{\rm d}t - \\ 
& \int_{\overset{r-1}{T}}^{\overset{r}{T}}\left|\zf\right|^2{\rm d}t-(\ln 2\pi-1-c)\overset{r-1}{T}\sim -c_0, \\ 
& r=1,\dots,k,\ \ln 2\pi -1-c\approx 0.26, 
\end{split}
\ee 
for every fixed $k\in\mbb{N}$, where 
\be \label{1.11} 
\overset{r-1}{T}=T^{r-1}(T)=\vp_1^{1-r}(T), 
\ee  
$c$ is the Euler's constant and $c_0$ is the Titchmarsh-Kober-Atkinson constant.  

\begin{remark}
An asymptotic \emph{conservation law} is expressed by the formula (\ref{1.10}). This conservation law concerns the behaviour of corresponding areas at $T\to\infty$, where 
\be \label{1.12} 
\overset{r-1}{T}=\int_0^{\overset{r-1}{T}} 1 {\rm d}t, 
\ee  
of course. This also gives the geometric meaning to the constant $c_0$. 
\end{remark} 

\begin{remark}
The area 
\be \label{1.13} 
\int_{\overset{r-1}{T}}^{\overset{r}{T}}\left|\zf\right|^2{\rm d}t 
\ee  
may be viewed as \emph{an external creator} of the asymptotic equilibrium (\ref{1.9}) for the complementary subsystems of three areas on the interval $(0,\overset{r-1}{T}]$. 
\end{remark} 

\begin{remark}
The area 
\be \label{1.14} 
\int_{t\in (0,\overset{r-1}{T}]^-}\left\{ -\ln t+\left|\zf\right|^2\right\}{\rm d}t
\ee 
is connected with the Titchmarsh $\Omega$-theorem in that sense, that it gives information about \emph{the outlets} of the values of $|\zf|^2$ above the curve $y=\ln t$. 
\end{remark}

\subsection{} 

In this paper we use the following notions of our works \cite{2} -- \cite{5}: 
\begin{itemize}
	\item[{\tt (a)}] Jacob's ladder $\vp_1(T)$, 
	\item[{\tt (b)}] direct iterations of Jacob's ladders 
	\bdis 
	\begin{split}
		& \vp_1^0(t)=t,\ \vp_1^1(t)=\vp_1(t),\ \vp_1^2(t)=\vp_1(\vp_1(t)),\dots , \\ 
		& \vp_1^k(t)=\vp_1(\vp_1^{k-1}(t))
	\end{split}
	\edis 
	for every fixed natural number $k$, 
	\item[{\tt (c)}] reverse iterations of Jacob's ladders 
	\bdis  
	\begin{split}
		& \vp_1^{-1}(T)=\overset{1}{T},\ \vp_1^{-2}(T)=\vp_1^{-1}(\overset{1}{T})=\overset{2}{T},\dots, \\ 
		& \vp_1^{-r}(T)=\vp_1^{-1}(\overset{r-1}{T})=\overset{r}{T},\ r=1,\dots,k, 
	\end{split} 
	\edis  
	where, for example, 
	\be \label{1.15} 
	\vp_1(\overset{r}{T})=\overset{r-1}{T}
	\ee  
	for every fixed $k\in\mbb{N}$, where 
	\be \label{1.16}
	\begin{split} 
		& \overset{r}{T}-\overset{r-1}{T}\sim(1-c)\pi(\overset{r}{T});\ \pi(\overset{r}{T})\sim\frac{\overset{r}{T}}{\ln \overset{r}{T}},\ r=1,\dots,k,\ T\to\infty, \\ 
		& \overset{0}{T}=T<\overset{1}{T}(T)<\overset{2}{T}(T)<\dots<\overset{k}{T}(T), \\ 
		& T\sim \overset{1}{T}\sim \overset{2}{T}\sim \dots\sim \overset{k}{T},\ T\to\infty,  
	\end{split}
	\ee 
	and $\pi(x)$ stands for the prime-counting function. 
\end{itemize} 

\begin{remark}
	The asymptotic behaviour of the points 
	\bdis 
	\{T,\overset{1}{T},\dots,\overset{k}{T}\}
	\edis  
	is as follows: at $T\to\infty$ these points recede unboundedly each from other and all together are receding to infinity. Hence, the set of these points behaves at $T\to\infty$ as one-dimensional Friedmann-Hubble expanding Universe. 
\end{remark}  

\section{The eleventh point of contact with the Fermat-Wiles theorem} 

\subsection{} 

The substitution 
\be \label{2.1} 
T\to \overset{r}{T},\ r=1,\dots,k 
\ee  
used in (\ref{1.9}), where\footnote{See (\ref{1.15}).} 
\bdis 
\vp_1(\overset{r}{T})=\overset{r-1}{T},\ \overset{r}{T}\sim T,\ T\to\infty
\edis 
gives us the following formula 
\be \label{2.2} 
\begin{split}
& \int_0^{\overset{r}{T}}\left|\zf\right|^2{\rm d}t=\overset{r-1}{T}\ln \overset{r-1}{T}+(c-\ln 2\pi)\overset{r-1}{T}+ \\ 
& c_0+\mcal{O}\left(\frac{\ln T}{T}\right). 
\end{split}
\ee 
Since 
\be \label{2.3} 
\overset{r-1}{T}\ln \overset{r-1}{T}=\int_0^{\overset{r-1}{T}}\ln t{\rm d}t+\overset{r-1}{T}, 
\ee  
we obtain the following Lemma. 

\begin{mydef51}
\be \label{2.4} 
\begin{split}
& \int_0^{\overset{r-1}{T}}\left\{ \ln t-\left|\zf\right|^2\right\}{\rm d}t-\int_{\overset{r-1}{T}}^{\overset{r}{T}}\left|\zf\right|^2{\rm d}t= \\ 
& (\ln 2\pi-1-c)\overset{r-1}{T}-c_0+\mcal{O}\left(\frac{\ln T}{T}\right), \\ 
& r=1,\dots,k,\ T>T_0>0, 
\end{split}
\ee 
for every fixed $k\in\mbb{N}$ and every sufficiently big $T_0$. 
\end{mydef51} 

\subsection{} 

Let us remind our almost linear formula\footnote{See \cite{6}, (3.4), (3.6); $0<\delta$ being sufficiently small.} 
\be \label{2.5} 
\int_{\overset{r-1}{T}}^{\overset{r}{T}}\left|\zf\right|^2{\rm d}t=(1-c)\overset{r-1}{T}+\mcal{O}(T^{1/3+\delta})
\ee 
for every fixed $k\in\mbb{N}$. 

\begin{remark}
Certain continuum set of increments of the Hardy-Littlewood integral (\ref{1.1}) is expressed by our formula (\ref{2.5}). These increments are constructed on the following set of segments\footnote{Comp. (\ref{1.16}).} 
\be \label{2.6} 
\begin{split}
& [\overset{r-1}{T}(T),\overset{r}{T}(T)]:\ \overset{r}{T}(T)-\overset{r-1}{T}(T)\sim (1-c)\pi(\overset{r}{T}), \\ 
& \overset{r}{T}(T)=\vp_1^{-r}(T),\ r=1,\dots,k,\ \overset{0}{T}=T,\ T\to\infty. 
\end{split}
\ee 
\end{remark} 

Now, using (\ref{2.5}) in (\ref{2.4}), we obtain the following lemma. 

\begin{mydef52}
\be \label{2.7} 
\begin{split}
& \int_0^{\overset{r-1}{T}}\left\{ \ln t-\left|\zf\right|^2\right\}{\rm d}t= \\ 
& = (\ln 2\pi-2c)\overset{r-1}{T}+\mcal{O}(T^{1/3+\delta}),\ r=1,\dots,k, 
\end{split}
\ee 
where $\ln 2\pi-2c\approx 0.68$. 
\end{mydef52} 

\subsection{} 

Next, we use (\ref{2.7}) with $r=1$, that is 
\be \label{2.8} 
\begin{split}
& \int_0^{T}\left\{ \ln t-\left|\zf\right|^2\right\}{\rm d}t= \\ 
& = (\ln 2\pi-2c)T+\mcal{O}(T^{1/3+\delta}),\ T>T_0. 
\end{split}
\ee 
If we substitute 
\be \label{2.9} 
T=\frac{x}{\ln 2\pi -2c}\tau;\ \tau>\frac{\ln2\pi-2c}{x}T_0\equiv \tau_0(x)
\ee 
into eq. (\ref{2.8}) for every fixed positive $x$, we obtain 
\be \label{2.10} 
\begin{split}
& \frac{1}{\tau}\int_0^{\frac{x}{\ln 2\pi -2c}\tau}\left\{ \ln t-\left|\zf\right|^2\right\}{\rm d}t= \\ 
& x+\mcal{O}_x(\tau^{-2/3+\delta}), 
\end{split}
\ee 
where the constant in the $\mcal{O}_x$-estimate depends on $x$. 

Consequently the following lemma holds true. 

\begin{mydef53}
\be \label{2.11} 
\lim_{\tau\to\infty}\frac{1}{\tau}\int_0^{\frac{x}{\ln 2\pi -2c}\tau}\left\{ \ln t-\left|\zf\right|^2\right\}{\rm d}t=x 
\ee 
for every fixed positive $x$. 
\end{mydef53} 

\subsection{} 

From (\ref{2.11}) we have 
\be \label{2.12} 
\begin{split}
& \lim_{\tau\to\infty}\frac{1}{\tau}\int_0^{\frac{x^n+y^n}{z^n}\frac{\tau}{\ln 2\pi-2c}}\left\{ \ln t- \left|\zf\right|^2 \right\}{\rm d}t = \\ 
& \frac{x^n+y^n}{z^n}. 
\end{split}
\ee 
Therefore the following theorem holds true. 

\begin{mydef11}
It is true that the condition 
\be \label{2.13} 
\begin{split}
\lim_{\tau\to\infty}\frac{1}{\tau}\int_0^{\frac{x^n+y^n}{z^n}\frac{\tau}{\ln 2\pi-2c}}\left\{ \ln t- \left|\zf\right|^2 \right\}{\rm d}t \not=1
\end{split}
\ee 
on all Fermat's rationals 
\bdis 
\frac{x^n+y^n}{z^n};\ x,y,z\in\mbb{N},\ n\geq 3
\edis 
represents the logarithmically modified $\zeta$-equivalent of the Fermat-Wiles theorem. 
\end{mydef11} 

\begin{remark}
We have used the formula (\ref{2.8}) as the special case of (\ref{2.7}) ($r=1$). Of course, we may have used also other cases $r=2,\dots,k$ by means of the substitution 
\be \label{2.14} 
\overset{r-1}{T}(T)=\rho;\ T=\vp_1^{r-1}(\rho), 
\ee  
where 
\be \label{2.15} 
T\to\infty \ \Leftrightarrow \ \rho\to\infty,  
\ee 
applied into (\ref{2.7}). This procedure leads finally to the formula 
\be \label{2.16} 
\begin{split}
& \int_0^\rho\left\{ \ln t- \left|\zf\right|^2 \right\}{\rm d}t = \\ 
& (\ln2\pi - 2c)\rho+\mcal{O}(\rho^{1/3+\delta}). 
\end{split}
\ee 
\end{remark} 

\section{Conservation laws for areas connected with the oscillations of the logarithmically modified function $|\zf|^2$} 

\subsection{} 

Next statement follows immediately from our Lemma 1 and eq. (\ref{1.7}). 

\begin{mydef54}
\be \label{3.1} 
\begin{split}
& \int_{t\in (0,\overset{r-1}{T}]^+}\left\{ \ln t-\left|\zf\right|^2\right\}{\rm d}t=\\ 
& \int_{t\in (0,\overset{r-1}{T}]^-}\left\{ -\ln t+\left|\zf\right|^2\right\}{\rm d}t+\int_{\overset{r-1}{T}}^{\overset{r}{T}}\left|\zf\right|^2{\rm d}t+\\ 
& (\ln2\pi-1-c)\overset{r-1}{T}-c_0+\mcal{O}\left(\frac{\ln T}{T}\right),\ T\to\infty. 
\end{split}
\ee 
\end{mydef54} 

And consequently, by means of the formula (\ref{3.1}), we obtain the following theorem. 

\begin{mydef12} 
\be \label{3.2} 
\begin{split}
& \lim_{T\to\infty}\left[
\int_{t\in (0,\overset{r-1}{T}]^+}\left\{ \ln t-\left|\zf\right|^2\right\}{\rm d}t-\right. \\ 
& \left. \int_{t\in (0,\overset{r-1}{T}]^-}\left\{ -\ln t+\left|\zf\right|^2\right\}{\rm d}t - \right. \\ 
& \left. \int_{\overset{r-1}{T}}^{\overset{r}{T}}\left|\zf\right|^2{\rm d}t-(\ln2\pi-1-c)\overset{r-1}{T}
\right]=-c_0, \\ 
& r=1,\dots,k,\ \ln2\pi-1-c\approx 0.26. 
\end{split}
\ee 
\end{mydef12}  

\begin{remark}
The limit law for behaviour of corresponding areas is expressed by the eq. (\ref{3.2}). Simultaneously, an exact expression for the constant $c_0$ is given by (\ref{3.2}). 
\end{remark}

\subsection{} 

The translation $r\to r+1$ applied in (\ref{2.4}) gives us the following formula 
\be \label{3.3} 
\begin{split}
& \int_{0}^{\overset{r}{T}}\left\{ \ln t-\left|\zf\right|^2\right\}{\rm d}t-\int_{\overset{r}{T}}^{\overset{r+1}{T}}\left|\zf\right|^2{\rm d}t= \\ 
& (\ln2\pi-1-c)\overset{r}{T}-c_0+\mcal{O}\left(\frac{\ln T}{T}\right). 
\end{split}
\ee 
Now, subtraction of (\ref{2.4}) from (\ref{3.3}) implies the following statement. 

\begin{mydef55}
\be \label{3.4} 
\begin{split}
& \int_{\overset{r-1}{T}}^{\overset{r}{T}}\left\{ \ln t-\left|\zf\right|^2\right\}{\rm d}t- \\ 
& \int_{\overset{r}{T}}^{\overset{r+1}{T}}\left|\zf\right|^2{\rm d}t+\int_{\overset{r-1}{T}}^{\overset{r}{T}}\left|\zf\right|^2{\rm d}t= \\ 
& (\ln2\pi-1-c)(\overset{r}{T}-\overset{r-1}{T})+\mcal{O}\left(\frac{\ln T}{T}\right), 
\end{split}
\ee 
and, of course, (comp. (\ref{1.7}) and (\ref{3.1})) 
\be \label{3.5} 
\begin{split}
& \int_{t\in(\overset{r-1}{T},\overset{r}{T}]^+}\left\{ \ln t-\left|\zf\right|^2\right\}{\rm d}t+\int_{\overset{r-1}{T}}^{\overset{r}{T}}\left|\zf\right|^2{\rm d}t=\\ 
& \int_{t\in(\overset{r-1}{T},\overset{r}{T}]^-}\left\{ -\ln t+\left|\zf\right|^2\right\}{\rm d}t+\int_{\overset{r}{T}}^{\overset{r+1}{T}}\left|\zf\right|^2{\rm d}t+ \\ 
& (\ln2\pi-1-c)(\overset{r}{T}-\overset{r-1}{T})+\mcal{O}\left(\frac{\ln T}{T}\right), \\ 
& r=1,\dots,k,\ T\to\infty. 
\end{split}
\ee 
\end{mydef55} 

And finally, next theorem follows from (\ref{3.5}).  

\begin{mydef13}
\be \label{3.6} 
\begin{split}
	& \lim_{T\to\infty}\left[
	\int_{t\in (\overset{r-1}{T},\overset{r}{T}]^+}\left\{ \ln t-\left|\zf\right|^2\right\}{\rm d}t-\right. \\ 
	& \left. \int_{t\in (\overset{r-1}{T},\overset{r}{T}]^-}\left\{ -\ln t+\left|\zf\right|^2\right\}{\rm d}t + \right. \\ 
	& \left. \int_{\overset{r-1}{T}}^{\overset{r}{T}}\left|\zf\right|^2{\rm d}t-\int_{\overset{r}{T}}^{\overset{r+1}{T}}\left|\zf\right|^2{\rm d}t-\right. \\ 
	& \left. (\ln2\pi-1-c)(\overset{r}{T}-\overset{r-1}{T})
	\right]=0, \\ 
	& r=1,\dots,k. 
\end{split}
\ee 
\end{mydef13} 

\begin{remark}
The \emph{0-limit law}  for the behaviour of corresponding areas at $T\to\infty$ is expressed by the eq. (\ref{3.6}). 
\end{remark} 

\section{The twelfth point of contact with the Fermat-Wiles theorem} 

\subsection{} 

Let us remind\footnote{See \cite{7}, (4.6).} the formula 
\be \label{4.1} 
\lim_{\tau\to\infty}\frac{1}{\tau}\int_{\frac{x}{1-c}\tau}^{[\frac{x}{1-c}\tau]^1}\left|\zf\right|^2{\rm d}t=x,\ x>0. 
\ee  
The quotient of eqs. (\ref{2.11}) and (\ref{4.1})\footnote{Factors are cancelled, of course.} gives the following result 
\be\label{4.2} 
\begin{split}
& \int_0^{\frac{x}{\ln2\pi-2c}\tau}\left\{ \ln t-\left|\zf\right|^2\right\}{\rm d}t\sim \\ 
& \int_{\frac{x}{1-c}\tau}^{[\frac{x}{1-c}\tau]^1}\left|\zf\right|^2{\rm d}t,\ \tau\to\infty 
\end{split}
\ee 
for every fixed positive $x$. 

\subsection{} 

Next, the subtract of eq. (\ref{4.2}) in the case $x\to mx, (m+1)x$, gives following set of formulas 
\be \label{4.3} 
\begin{split}
& \int_{\frac{mx}{\ln2\pi-2c}\tau}^{\frac{(m+1)x}{\ln2\pi-2c}\tau}\left\{ \ln t-\left|\zf\right|^2\right\}{\rm d}t	\sim \\ 
& \int_{\frac{(m+1)x}{1-c}\tau}^{[\frac{(m+1)x}{1-c}\tau]^1}\left|\zf\right|^2{\rm d}t-\int_{\frac{mx}{1-c}\tau}^{[\frac{mx}{1-c}\tau]^1}\left|\zf\right|^2{\rm d}t, \\ 
& m=1,\dots,N,\ \tau\to\infty 
\end{split}
\ee 
for every fixed $N\in\mbb{N}$. 

Summing the system of formulas (\ref{4.3}) over $m$ in given range results in the following: 
\be \label{4.4} 
\begin{split}
& \int_{\frac{x}{\ln2\pi-2c}\tau}^{\frac{(N+1)x}{\ln2\pi-2c}\tau}\left\{ \ln t-\left|\zf\right|^2\right\}{\rm d}t\sim \\ 
& \int_{\frac{(N+1)x}{1-c}\tau}^{[\frac{(N+1)x}{1-c}\tau]^1}\left|\zf\right|^2{\rm d}t - 
\int_{\frac{x}{1-c}\tau}^{[\frac{x}{1-c}\tau]^1}\left|\zf\right|^2{\rm d}t
\end{split}
\ee 
for every fixed $N\in\mbb{N}$ and $\tau\to\infty$. Since it is true that 
\be \label{4.5} 
\begin{split}
& \int_{\frac{(N+1)x}{1-c}\tau}^{[\frac{(N+1)x}{1-c}\tau]^1}\left|\zf\right|^2{\rm d}t\sim (N+1)x\tau, \\ 
& \int_{\frac{x}{1-c}\tau}^{[\frac{x}{1-c}\tau]^1}\left|\zf\right|^2{\rm d}t\sim x\tau, 
\end{split}
\ee 
then, by the almost linear formula (\ref{2.5}), we obtain the following statement. 

\begin{mydef56}
\be \label{4.6} 
\frac{1}{\tau}\int_{\frac{x}{\ln2\pi-2c}\tau}^{\frac{(N+1)x}{\ln2\pi-2c}\tau}\left\{ \ln t-\left|\zf\right|^2\right\}{\rm d}t\sim Nx,\ \tau\to\infty 
\ee 
for every fixed positive $x$ and every fixed $N\in\mbb{N}$. 
\end{mydef56} 

\subsection{} 

From (\ref{4.6}) it follows immediately: 

\begin{mydefII1}
It is true that 
\be \label{4.7} 
\begin{split}
& \lim_{\tau\to\infty}\frac{1}{\tau}\int_{\frac{x^n+y^n}{z^n}\frac{\tau}{\ln2\pi-2c}}^{(N+1)\frac{x^n+y^n}{z^n}\frac{\tau}{\ln2\pi-2c}}\left\{ \ln t-\left|\zf\right|^2\right\}{\rm d}t= N\frac{x^n+y^n}{z^n}
\end{split}
\ee 
for every fixed $N\in\mbb{N}$ and every fixed Fermat's rational 
\bdis 
\frac{x^n+y^n}{z^n}. 
\edis 
\end{mydefII1} 

Finally, the following theorem holds true. 

\begin{mydef14}
The condition 
\be \label{4.8} 
\lim_{\tau\to\infty}\frac{1}{\tau}\int_{\frac{x^n+y^n}{z^n}\frac{\tau}{\ln2\pi-2c}}^{(N+1)\frac{x^n+y^n}{z^n}\frac{\tau}{\ln2\pi-2c}}\left\{ \ln t-\left|\zf\right|^2\right\}{\rm d}t\not=N 
\ee  
on all Fermat's rationals 
\bdis 
\frac{x^n+y^n}{z^n},\ x,y,z,n\in\mbb{N},\ n\geq 3
\edis 
separately for every fixed $N\in\mbb{N}$ represents the second variant of the logarithmically modified $\zeta$-equivalent of the Fermat-Wiles theorem. 
\end{mydef14} 

\section{On localized logarithmically modified $\zeta$-equivalent of the Fermat-Wiles theorem} 

\subsection{} 

If we use the section 7 of our work \cite{10} then we obtain next corollary to the Theorem 1. 

\begin{mydef41} 
It is true that the condition 
\be \label{5.1} 
\lim_{\tau\to\infty}\frac{1}{\tau}\int_{0}^{\frac{x^n+y^n}{z^n}\frac{\tau}{\ln2\pi-2c}}\left\{ \ln t-\left|\zf\right|^2\right\}{\rm d}t\not=1
\ee 
for 
\be \label{5.2} 
\forall\- \frac{x^n+y^n}{z^n}\in (1-\epsilon,1+\epsilon)
\ee 
represents the localized logarithmically modified $\zeta$-equivalent of the Fermat-Wiles theorem for every fixed and small $\epsilon>0$. 
\end{mydef41}

And once again, Theorem 4 implies 

\begin{mydef42}
It is true that the condition 
\be \label{5.3} 
\lim_{\tau\to\infty}\frac{1}{\tau}\int_{\frac{x^n+y^n}{z^n}\frac{\tau}{\ln2\pi-2c}}^{(N+1)\frac{x^n+y^n}{z^n}\frac{\tau}{\ln2\pi-2c}}\left\{ \ln t-\left|\zf\right|^2\right\}{\rm d}t\not=N
\ee 
for 
\be \label{5.4} 
\forall\- \frac{x^n+y^n}{z^n}\in \left(1-\frac{\epsilon}{N},1+\frac{\epsilon}{N}\right)
\ee 
represents the localized logarithmically modified $\zeta$-equivalent of the Fermat-Wiles theorem for every fixed $N\in\mbb{N}$. 
\end{mydef42}

I would like to thank Michal Demetrian for his moral support of my study of Jacob's ladders.

\end{document}